\documentclass{amsart}[12pt]

\newcommand{\N}{{\mathbb{N}}}

\newcommand{\R}{{\mathbb{R}}}

\newcommand{\Bh}{{\mathcal B}}
\newcommand{\Ch}{{\mathcal C}}
\newcommand{\Dh}{{\mathcal D}}

\newcommand{\Zh}{{\mathcal Z}}

\newcommand{\be}{\mathbf{1}}

\newcommand{\dr}{\mathrm{dr}\,}

\newcommand{\Prim}{\mathrm{Prim}}

\newcounter{number}[section]

\newenvironment{nummer}{\refstepcounter{number}{\noindent\arabic{section}.\arabic{number}}}{}

\newcommand{\bn}{\noindent \begin{nummer} \rm}
\newcommand{\en}{\end{nummer}}

\newenvironment{thms}{\noindent {\sc Theorem:} \it}{}

\newenvironment{props}{\noindent {\sc Proposition:} \it}{}
\newenvironment{dfs}{\noindent {\sc Definition:} \it}{}
\newenvironment{cors}{\noindent {\sc Corollary:} \it}{}

\newenvironment{nproof}{\noindent {\sc Proof:}}{\mbox{}\hfill 
\rule[-.2ex]{.25em}{1.8ex}}

\begin{document}

\title[A note on subhomogeneous $C^{*}$-algebras]{A note on subhomogeneous $C^{*}$-algebras}

\author{Ping Wong Ng}
\address{The Fields Institute for Research in Mathematical Sciences\\
222 College Street\\ \indent Toronto, Ontario\\ M5T 3J1\\ Canada}

\email{pwn@erdos.math.unb.ca}

\author{Wilhelm Winter}
\address{Mathematisches Institut der Universit\"at M\"unster\\
Einsteinstr. 62\\ D-48149 M\"unster}

\email{wwinter@math.uni-muenster.de}

\date{December 21, 2005}
\thanks{{\it Supported by:} DFG (through the SFB 478), EU-Network  Quantum Spaces - Noncommutative 
\indent Geometry (Contract No. HPRN-CT-2002-00280)}

\setcounter{section}{-1}

\begin{abstract}
We show that finitely generated subhomogeneous $C^{*}$-algebras have finite decomposition rank. As a consequence,  any separable ASH $C^{*}$-algebra can be written as an inductive limit of subhomogeneous $C^{*}$-algebras each of which has finite decomposition rank. \\
It then follows from work of H.\ Lin and of the second named author that the class of simple unital ASH algebras which have real rank zero and absorb the Jiang--Su algebra tensorially satisfies the Elliott conjecture.
\end{abstract}

\maketitle

\section{Introduction}

\noindent
Any compact metrizable space $X$ can be approximated by topologically finite-dimensional spaces in the sense that it can be written as an inverse limit of  polyhedra $S_{n}$, each of which is  finite-dimensional (in general, the dimension of the $S_{n}$ will not be globally bounded, though). For the corresponding $C^{*}$-algebras this means that $\Ch(X)$ is the inductive limit of  $C^{*}$-algebras $\Ch(S_{n})$ each of which is topologically finite-dimensional, and the respective statement holds for   $r\times r$-matrices over $\Ch(X)$. It is natural to ask whether this construction can be carried over to more general $C^{*}$-algebras. \\
For example, it is not hard to show that if $A$ is separable and homogeneous (i.e., all irreducible representations of $A$ have the same rank), then one can write $A$ as an inductive limit of homogeneous $C^{*}$-algebras $A_{n}$ each of which has topologically finite-dimensional spectrum. In \cite{Bl2}, Proposition 2.2, Blackadar has outlined a proof of a more general fact: if $A$ is separable and approximately homogeneous (AH), i.e., an inductive limit of  homogeneous $C^{*}$-algebras, then $A$ can be written as an inductive limit of homogeneous $C^{*}$-algebras of finite topological dimension. In these notes we shall give a simple proof of the respective statements for (approximately)  subhomogeneous $C^{*}$-algebras. Recall that a $C^{*}$-algebra is subhomogeneous if there is an upper bound for the ranks of its irreducible  representations. Inductive limits of such algebras are then called approximately subhomogeneous (ASH).  These results will be immediate consequences of our Theorem \ref{main-result} below, which says that finitely generated subhomogeneous $C^{*}$-algebras have finite topological dimension.\\
In general it is not at all clear how to define the  topological dimension of a $C^{*}$-algebra (there are several choices), but for a homogeneous $C^{*}$-algebra $A$ it seems natural to  define the topological dimension as covering dimension of $\Prim A$ (the space of kernels of irreducible representations), $\dim \Prim A$ (see \cite{HW} for a general account on covering dimension). Note that, for homogeneous $A$, $\Prim A$ coincides with the spectrum $\hat{A}$, so there is no ambiguity in the choice of the underlying space. This definition can easily be generalized to the subhomogeneous case by defining the topological dimension to be $\max_{k \in \N} \{\dim  \Prim_{k} A\}$, where $\Prim_{k}A$ is the space of kernels of irreducible representations of rank $k$. \\
For (sub)homogeneous $C^{*}$-algebras, the  topological dimension coincides with the decomposition rank of $A$, $\dr A$, a notion of covering dimension which makes sense for arbitrary stably finite nuclear $C^{*}$-algebras (cf.\ \cite{KW}, \cite{Wi1}  and \cite{Wi3}). In \cite{Wi6}, the second named author has introduced the notion of locally finite decomposition rank. It is clear from \cite{Bl2}, Proposition 2.2, that AH algebras have locally finite decomposition rank, and it will follow from the results in the present paper that this also holds for ASH algebras. In this sense, locally finite decomposition rank generalizes the concepts of AH and ASH algebras.\\
Our main source of interest in approximating subhomogeneous $C^{*}$-algebras by topologically finite-dimensional ones stems from the relevance of this problem for Elliott's classification program. This is an attempt to classify nuclear $C^{*}$-algebras by $\mathrm{K}$-theoretic invariants. At least for AH algebras this often leads to the problem of comparing vector bundles over compact spaces. Although this is not possible in general, nonstable $\mathrm{K}$-theory provides partial answers in the case of base spaces with finite covering dimension. For simple AH algebras with globally bounded topological dimension one then obtains good comparison properties, and  there are similar results for $C^{*}$-algebras with finite decomposition rank (cf.\ \cite{Wi4}). Building on the ideas of \cite{Wi5}, \cite{Wi6} provides a classification theorem for simple, unital, separable $C^{*}$-algebras with real rank zero and locally finite decomposition rank which also absorb the Jiang--Su algebra $\Zh$ tensorially. Together with our result this means in particular that simple, unital, separable ASH algebras with real rank zero are classified by the Elliott invariant up to $\Zh$-stability.

\section{The main result}

\bn
\begin{dfs}
We say a $C^{*}$-algebra $A$ is $r$-subhomogeneous, if there is $r \in \N$ such that every irreducible representation of $A$ has rank less than or equal to $r$. \\
We say $A$ is $r$-homogeneous, if every irreducible representation has rank equal to $r$. Sometimes we shall not specify the $r$ explicitly. 
\end{dfs}
\en

\bn
For a $C^{*}$-algebra $A$ let $\Prim A$ denote the primitive ideal space. If $X$ is a closed subset of $\Prim A$, then $J_{X}:=\bigcap_{x \in \Prim A \setminus X} x$ is a closed two-sided ideal of $A$; we denote the quotient $A/J_{X}$ by $A_{X}$.  The map $X \mapsto J_{X}$ is a bijection between the set of closed subsets of $\Prim A$ and the set of closed ideals of $A$.\\
For $k \in \N$, let $\Prim_{k} A \subset \Prim A$ be the subspace of kernels of $k$-dimensional irreducible representations. Recall from \cite{Dx}, Chapter 3,  that $\Prim_{k}A $ is a locally compact Hausdorff space. If $A$ is separable, the spaces $\Prim_{k} A$ are second countable. \\
If $A$ is $r$-subhomogeneous, then $\Prim A = \bigcup_{k=1}^{r} \Prim_{k}A$ and we may define the topological dimension of $A$  as $\max_{k=1, \ldots,r} \{\dim \Prim_{k} A \}$. Using results of \cite{P4}, it was shown in \cite{Wi3} that this number coincides with $\dr A$, the decomposition rank of $A$. In \cite{P4} it was shown that a subhomogeneous $C^{*}$-algebra with finite topological dimension can be written as an iterated pullback of homogeneous $C^{*}$-algebras with finite topological dimension (see \cite{P4} for a precise statement and for more information on these so-called recursive subhomogeneous $C^{*}$-algebras). 
\en

\bn
We shall have use for an easy consequence of \cite{Fe}, Theorem 3.2; it is well-known, but we include a proof for completeness. Essentially, it says that the restriction of $A$ to $\Prim_{k}A$, when regarded as a bundle over $\Prim_{k}A$ with fibres $M_{k}$, is locally trivial. 

\label{compact-neighborhood}
\begin{props}
Let $A$ be a separable $C^{*}$-algebra and $k \in \N$. Then, any $x \in \Prim_{k} A$ has a compact neighborhood (with respect to the relative topology of $\Prim_{k}A$) $X$ in $\Prim_{k} A$ such that $A_{X} \cong \Ch(X) \otimes M_{k}$. 
\end{props}

\begin{nproof}
Since $\Prim_{k}A$ is a locally compact Hausdorff space, $x$ has a compact neighborhood $Y \subset \Prim_{k}A$.  $Y$ clearly is a closed subset of $\Prim A$; moreover,  $A_{Y}$ is $k$-homogeneous: Indeed, if $\pi$ is an irreducible representation of $A_{Y}$, its composition with the quotient map yields an irreducible representation  (also denoted by $\pi$) of $A$ which annihilates all of $J_{Y}$, whence $J_{Y} \subset \ker \pi =:y \in \Prim A$. This implies $J_{Y \cup \{y\}} = (\bigcap_{x \in Y} x) \cap y = J_{Y} \cap y = J_{Y}$. Since the map $Y \mapsto J_{Y}$ is a bijection, this means that $Y \cup \{y\}= Y$; therefore, $y \in Y \subset \Prim_{k} A$; it follows that the representation $\pi$ is $k$-dimensional.   \\
Now \cite{Fe}, Theorem 3.2, says that $A_{Y}$ is isomorphic to the algebra $\Ch(\Bh)$ of continuous sections of a fibre bundle $\Bh$ with base space $Y$ and fibres $M_{k}$. Since $\Bh$ is locally trivial, there is a compact neighborhood $X \subset Y$ of $x$ such that the restriction of $\Bh$ to $X$ is of the form $\Bh_{X} \cong X \times M_{k}$. Note that the section algebra  $\Ch(\Bh_{X})$ is nothing but  $\Ch(X,M_{k}) \cong \Ch(X) \otimes M_{k}$.  On the other hand, $A_{X} = (A_{Y})_{X} \cong \Ch(\Bh)_{X} \cong \Ch(\Bh_{X})$, since the restriction map $\Ch(\Bh) \to \Ch(\Bh_{X})$ is surjective. 
\end{nproof}
\en

\bn
\label{finitely-generated-homogeneous}
\begin{props}
Let $X$ be a compact space and let $k,m \in \N$ be given. If $\Ch(X) \otimes M_{k}$ is  generated by $m$ elements, then $\dim X \le 4 \cdot m \cdot k^{2}$.
\end{props}

\begin{nproof}
Let  $a^{(1)}, \ldots, a^{(m)}$ be the $m$ generators of $\Ch(X) \otimes M_{k}$, and let $a^{(l)}_{ij} \in \Ch(X)$, $i,j =1, \ldots, k$, $l=1, \ldots, m$ denote the coordinate functions; since these separate points of $X$, together with $\be_{X}$ they generate $\Ch(X)$ as a $C^{*}$-algebra. By taking positive and negative parts of the real and imaginary parts of the $a^{(l)}_{ij}$ we obtain  positive functions $b_{j}$, $j=1, \ldots, n:=4\cdot m \cdot k^{2}$ which (together with $\be_{X}$) also generate $\Ch(X)$ as a $C^{*}$-algebra. By rescaling the $b_{j}$ we may clearly assume that $0 \le \sum_{j=1}^{n} b_{j} \le \be_{X}$ for all $j$. Now $b_{1}, \ldots, b_{n}, \be_{X} - \sum_{j=1}^{n} b_{j}$ are $n+1$ positive elements which add up to $\be_{X}$ and generate $\Ch(X)$ as a $C^{*}$-algebra. Let $\Delta^{n} \subset \R^{n+1}$ denote the standard simplex with $n+1$ vertices. Since $\Ch(\Delta^{n})$ is the universal unital $C^{*}$-algebra generated by $n+1$ pairwise commuting positive elements which add up to the unit,  we see that $\Ch(X)$ is a quotient of $\Ch(\Delta^{n})$, whence $X$ is homeomorphic to a closed subset of $\Delta^{n}$. But then by \cite{HW}, Theorem III.1, $\dim X \le \dim \Delta^{n} = n = 4\cdot m \cdot k^{2}.$
\end{nproof}
\en

\bn
\label{main-result}
\begin{thms}
Let $A$ be a finitely generated subhomogeneous $C^{*}$-algebra. Then, $\dr A < \infty$. If $A$ is unital, it is recursive subhomogeneous of finite topological dimension in the sense of \cite{P4}.
\end{thms}

\begin{nproof}
Suppose $A$ is generated by $m$ elements. For $k \in \N$, apply Proposition \ref{compact-neighborhood} to each $x \in \Prim_{k}A$ to obtain a compact neighborhood $X_{x} \subset \Prim_{k} A$ such that $A_{X_{x}} \cong \Ch(X_{x}) \otimes M_{k}$. Since $A_{X_{x}}$ is a quotient of $A$, it is also generated by $m$ elements and we conclude from Proposition \ref{finitely-generated-homogeneous} that $\dim X_{x} \le 4 \cdot m \cdot k^{2}$. We have that $\Prim_{k}A = \bigcup_{x \in \Prim_{k}A} X_{x}$. But since $A$ is finitely generated it is separable, whence $\Prim_{k}A$ is locally compact and second countable. Therefore, $\Prim_{k}A$ is covered by countably many of the $X_{x}$ and by the countable sum theorem for covering dimension (\cite{HW}, Theorem III.2) we have $\dim \Prim_{k} A \le \max_{x \in X} \dim X_{x} \le 4 \cdot m \cdot k^{2}$. Since $A$ is subhomogeneous, there is $r \in \N$ such that the $\Prim_{k}A$ are empty for $k > r$; it follows that $\max_{k} \{\dim \Prim_{k}A\}  \le 4 \cdot m \cdot r^{2} < \infty$. Now \cite{Wi3}, Theorem 1.6, says that $\dr A = \max_{k} \{\dim \Prim_{k}A\} < \infty$. Note that \cite {Wi3} uses Theorem 2.16 of \cite{P4} and that, in the unital case,  this result already implies that $A$ is recursive subhomogeneous of finite topological dimension.      
\end{nproof}
\en

\section{Some consequences}

\bn
\label{ASH-RSH}
\begin{cors}
Any separable ASH algebra $A$ may be written as an inductive limit of subhomogeneous $C^{*}$-algebras $A_{n}$ each of which has finite decomposition rank. If $A$ is unital, the $A_{n}$ may be chosen to be recursive subhomogeneous of finite topological dimension. 
\end{cors}

\begin{nproof}
Let $A = \lim_{\to} B_{n}$ be an inductive limit decomposition of $A$ such that each $A_{n}$ is subhomogeneous. Since quotients of subhomogeneous $C^{*}$-algebras are again subhomogeneous, we may replace the $B_{n}$ by their images in the limit, i.e., we may assume $A = \overline{\bigcup_{n \in \N} B_{n}}$ with $B_{n} \subset B_{n+1} \; \forall \, n \in \N$. Since $A$ is separable, it is straightforward to construct a sequence $(a_{n})_{n \in \N}$ which is dense in $A$ and satisfies $\{a_{0}, \ldots, a_{n}\} \subset B_{n} \; \forall \, n \in \N$. Set $A_{n}:= C^{*}(a_{0}, \ldots, a_{n}) \subset B_{n}$, then $A = \overline{\bigcup_{n \in \N} A_{n}}$ and $A_{n} \subset A_{n+1} \; \forall \, n \in \N$, i.e., $A$ is the inductive limit of the $A_{n}$. But each $A_{n}$ is subhomogeneous, as it is a subalgebra of a subhomogeneous algebra;  it is finitely generated, so it has finite decomposition rank by Theorem \ref{main-result}. If $A$ is unital, we may clearly assume each $A_{n}$ to be unital and the statement also follows from \ref{main-result}. 
\end{nproof}
\en

\bn
As an immediate consequence we have:

\begin{cors}
Any separable ASH algebra $A$ has locally finite decomposition rank in the sense of \cite{Wi6}.
\end{cors}
\en

\bn
Combining the preceding Corollary with the results of \cite{Wi6}  and \cite{Li1}, we obtain the following classification theorem (see \cite{R1} for an introduction to the Elliott program and for a description of the invariant):

\begin{cors}
The class of simple, separable, unital ASH algebras which have real rank zero and absorb the Jiang--Su algebra $\Zh$ tensorially satisfies the Elliott conjecture. In other words, two such algebras are isomorphic if and only if their Elliott invariants are,  and any isomorphism of the invariants is induced by an isomorphism of the algebras. 
\end{cors}
\en

\bn
We may also use Corollary \ref{ASH-RSH} to generalize \cite{TW1}, Corollaries 5.8 and 5.11:

\begin{cors}
Let $A$ be a separable ASH algebra which absorbs the Jiang--Su algebra $\Zh$ tensorially. Then, $A$ has strict slow dimension growth (in the sense of \cite{P3}). \\
If $\Dh$ is ASH and strongly self-absorbing in the sense of \cite{TW1}, then $\Dh$ is either UHF or projectionless. Moreover, $\Dh$ has strict slow dimension growth.
\end{cors}
\en

\end{document}